\documentclass [12pt]{article}

\usepackage{amssymb}
\usepackage{amsfonts}
\usepackage{graphicx, epsfig}
\usepackage{amsmath}
\usepackage{graphicx}
\usepackage{color}%
\providecommand{\U}[1]{\protect\rule{.1in}{.1in}}
\begin{document}

\begin{center}
\ \ \ {\Large Dynamic Principal Components in the Time Domain\vspace{1cm}}

Daniel Pe\~{n}a and Victor J. Yohai \vspace{7cm}
\end{center}

\hrule\smallskip\smallskip

\noindent Daniel Pe\~{n}a is Professor, Statistics Department, Universidad
Carlos III de Madrid, Calle Madrid 126, 28903 Getafe, Espa\~{n}a, (E-mail:
daniel.pe\~{n}a@uc3m.es). V\'{\i}ctor J. Yohai is Professor Emeritus,
Mathematics Department, Faculty of Exact Sciences, Ciudad Universitaria, 1428
Buenos Aires, Argentina (E-mail: victoryohai@gmail.com). This research was
partially supported by Grant ECO2012-38442 from MINECOM, Spain, and Grants
W276 from Universidad of Buenos Aires, PIP's 112-2008-01-00216 and
112-2011-01- 00339 from CONICET and PICT 2011-0397 from ANPCYT, Argentina.

\newpage

\noindent\ 

\begin{center}
ABSTRACT
\end{center}

We propose a time domain \ approach to define \ dynamic principal components
(DPC) using a reconstruction of the original series criterion. This approach
\ to define DPC was introduced by Brillinger, who gave a very elegant
theoretical solution in the stationary case using the cross spectrum. Our
procedure can be applied under more general conditions including the case of
non stationary series and relatively short series. We also present a robust
version of our procedure that allows to estimate the DPC when the series have
\ outlier contamination. \ \ Our non robust and robust procedures are
illustrated with real datasets.\ 


\section{Introduction}

Dimension reduction is very important in vector time series because the number
of parameters in a model grows very fast with the dimension $m$ of the vector
of time series. Therefore, finding simplifying structures or factors in these
models is important to reduce the number of parameters required to apply them
to real data. Besides, these factors, as we will see in this paper, may allow
to reconstruct with a small error the set of data and therefore reducing the
amount of information to be stored. In this article, we will consider linear
time series models and we will concentrate in the time domain approach.
Dimension reduction is usually achieved by finding linear combinations of the
time series variables which have interesting properties. Suppose the time
series vector $\mathbf{z}_{t}=(z_{1,t},...,z_{m,t})^{\prime}$, where $1\leq
t\leq T,$ and we assume, for simplicity, that $\overline{\mathbf{z}}%
=T^{-1}\sum_{t=1}^{T}\mathbf{z}_{t},$ which will estimate the mean if the
process is stationary, is zero. It is well known that the first principal
component, $p_{1,t},1\leq t\leq T,$ minimizes the mean squared prediction
error of the reconstruction of the vector time series, given by $\sum
_{j=1}^{m}\sum_{t=1}^{T}(z_{j,t}-\alpha_{j}p_{1,t})^{2}$ and, in general, the
first $k$ principal components, $k\leq m,$ $p_{1t},...,p_{kt},$ $1\leq t\leq
T,$ minimize the mean squared prediction error $\sum_{j=1}^{m}\sum_{t=1}%
^{T}(z_{j_{,}t}-\sum_{i=1}^{k}\alpha_{j,i}p_{i,t})^{2}$ to reconstruct the
vector of time series. Let $C=\sum_{t=1}^{T}\mathbf{z}_{t}\mathbf{z}%
_{t}^{\prime}/T,$ be the sample covariance matrix and let $\lambda_{1}%
\geq\lambda_{1}\geq\lambda_{m}$ be the eigenvalues \ of $C.$ Then
$\mathbf{\alpha}_{i}=(\alpha_{1,i},...,\alpha_{m,i})^{\prime},$ $1\leq i\leq
m,$ is the eigenvectors of $C$ corresponding to the eigenvalue $\lambda_{i}.$\ 

Ku, Storer and Georgakis (1995)\ propose to apply principal components to the
augmented observations $\mathbf{z}_{t}^{\ast}=(\mathbf{z}_{t-h}^{\prime
},\mathbf{z}_{t-h+1}^{\prime},...,\mathbf{z}_{t}^{\prime}$ )$^{\prime
}\mathbf{,}$ $h+1\leq t\leq T,\mathbf{\ }$\ that includes the values of the
series up to lag $h.$ These principal components provide linear combinations
of the present and past values of the time series with largest variance, and
using the well know properties of standard principal components we conclude
that the first component obtained from this approach is a solution to the
following reconstruction problem%
\[
M_{1}=\sum_{j=1}^{m}\left[  \sum_{t=h+1}^{T}(z_{j,t}-\alpha_{j}p_{t-h}%
)^{2}+\sum_{t=h}^{T-1}(z_{j,t}-\alpha_{j}p_{t-h+1})^{2}+...+\sum_{t=1}%
^{T-h}(z_{j,t}-\alpha_{j}p_{t})^{2}\right]  ,
\]
which implies that, apart from the end effect, we minimize for each
observation $z_{j,t}$, for $h+1\leq t\leq T-h,$ the sum $\sum_{j=1}^{m}%
\sum_{l=0}^{h}(z_{j,t}-\alpha_{j}p_{t-l})^{2}$ .Thus, this approach does not
optimize a useful reconstruction criterion.

An alternative way to find interesting linear combinations was proposed by Box
and Tiao (1977) who suggested maximizing the predictability of the linear
combinations $c_{t}=\mathbf{\gamma}^{\prime}\mathbf{z}_{t}$. \ Other linear
methods for dimension reduction in time series models have been given by the
scalar component models, SCM, (Tiao and Tsay, 1989), the reduced-rank models
(Ahn and Reinsel, 1990, Reinsel and Velu, 1998), and dynamic factor models
(Pe\~{n}a and Box, 1987, Stock and Watson, 1988, Forni el al. 2000, Pe\~{n}a
and Poncela 2006 and Lam and Yao 2012), among others. None of the previous
mentioned methods has as a goal to reconstruct the original series by using
the principal components as in the classical case.

Brillinger (1981) addressed the reconstruction problem as follows. Suppose now
the zero mean $m$ dimensional stationary process $\left\{  \mathbf{z}%
_{t}\right\}  ,$ $-\infty<t<\infty.$ Then, \ the dynamic principal components
\ are defined by searching for $m\times1$ vectors $\mathbf{c}_{h}%
,-\infty<h<\infty$ and $\mathbf{\beta}_{j},-\infty<j<\infty,$ so that if we
consider as first principal component the linear combination
\begin{equation}
f_{t}=%
{\displaystyle\sum\limits_{h=-\infty}^{\infty}}
\mathbf{c}_{h}^{\prime}\mathbf{z}_{t-h}, \label{BR1}%
\end{equation}
then%
\begin{equation}
E\left[  (\mathbf{z}_{t}-\sum_{j=-\infty}^{\infty}\mathbf{\beta}_{j}%
f_{t+j})^{\prime}(\mathbf{z}_{t}-\sum_{j=-\infty}^{\infty}\mathbf{\beta}%
_{j}f_{t+j})\right]  . \label{BR2}%
\end{equation}
is minimum. Brillinger elegantly solved this problem \ by showing that
$\mathbf{c}_{k}$ is the inverse Fourier transform of the principal components
of the \ cross spectral matrices for each frequency, and $\mathbf{\beta}_{j}$
is the inverse Fourier transform of the conjugates of the same principal
components. See Brillinger (1981) and Shumway and Stoffer (2000) for the
details of the method. Although this result solves the theoretical problem it
has the following shortcomings: (i) It can be applied only to stationary
series;\ (ii) The \ optimal solution requires the unrealistic assumption that
infinite series are observed, and it is not clear how to modify it\ when the
observed series are finite; (iii) It is not clear how to robustify these
principal components using a reconstruction criterion. The second shortcoming
seems specially serious. In fact in Section \ref{simul} we show by means of a
Monte Carlo simulation that what seems a natural \ modification for finite
series of the Brillinger's procedure does not work well.

In this paper we address the sample reconstruction of a vector of time series
avoiding the drawbacks of Brillinger method. Our procedure provides an optimal
reconstruction of the vector of time series from a finite number of lags. Some
of the advantages of our procedure are: (i) it does not require stationarity
and (ii) it can be easily made robust by changing the minimization of the mean
squared error criterion by the minimization of a robust scale. The rest of
this article is organized as follows. In Section \ref{recons} we describe the
proposed dynamic principal components based on the reconstruction criterion.
In Section \ref{RECOn1} we study the particular case where the proposed
dynamic principal components depend only on one lag. In Section \ref{simul} we
show the results of a Monte Carlo study that compares the proposed dynamic
principal components, with the ordinary principal components and those
proposed by Brillinger and we show the performances of these three types of
principal components in two\ real examples. In Section \ref{SecRDC} we define
robust dynamic principal components using a robust reconstruction criterion
and illustrate in one example the good performance of this estimator to
eliminate the influence of outliers. In Section \ref{section: con} some final
conclusions are presented. Section \ref{Sec: Ap} is an Appendix containing
mathematical derivations.

\section{Finding time series with optimal reconstruction properties
\label{recons}}

Suppose that \ we observe $z_{j,t},1\leq$ $j\leq m,$ $1\leq t\leq T,$ \ and
consider \ \ two integer numbers $k_{1}\geq0$ and $k_{2}\geq0.$ We can define
the first dynamic principal component with $\ k$ lags (first DPC$_{k}$) as a
vector $\mathbf{f=}(f_{t})_{-k_{1}+1\leq t\leq T+k_{2}},$ so \ that the
reconstruction of series $z_{j,t},1\leq$ $j\leq m,$ as a linear combination
\ of $f_{t-k_{1}},f_{t-k_{1}+1},....f_{t},f_{t+1},...,f_{t+k_{2}}$ is optimal
with the mean squared error (MSE) criterion. More precisely, suppose that
given a possible factor $\mathbf{f,}$\textbf{\ }the $m\times(k_{1}+k_{2})$
matrix of coefficients $\mathbf{\beta=(}\beta_{j,i})_{1\leq j\leq
m,-k_{1}+1\leq i\leq k_{2}},$and $\mathbf{\alpha=(}\alpha_{1},...,\alpha_{m})$
are used to reconstruct the values $z_{j,t}$ as
\[
\widehat{z}_{j,t}(\mathbf{f,\beta}_{j},\alpha_{j}\mathbf{)}=\sum_{i=-k_{1}%
}^{k_{2}}\beta_{j,i}f_{t+i}+\alpha_{j},
\]
where $\mathbf{\beta}_{j}$ is the $j$-th row of $\mathbf{\beta.}$ Let
$\ k=k_{1}+k_{2}$ and put $f_{t}^{\ast}=f_{t-k_{1}},$ $1\leq t\leq T+k$ and
$\beta_{j.,h}^{\ast}=$ $\beta_{j,h-k_{1}},0\leq h\leq k,$ then, the
reconstructed series are obtained as
\[
\widehat{z}_{j,t}(\mathbf{f,\beta}_{j},\alpha_{j}\mathbf{)}\mathbf{=}%
\sum_{i=-k_{1}}^{k}\beta_{j,i}f_{t+i+k_{1}}^{\ast}+\alpha_{j}=\sum_{h=0}%
^{k}\beta_{j,h,}^{\ast}f_{t+h}^{\ast}+\alpha_{j}.
\]
Therefore we can always assume that $k_{1}=0$ and we will use $\ k$ to denote
the number of forward lags.

Consider the MSE loss function
\begin{equation}
\text{MSE}(\mathbf{f,\beta,\alpha)}=\sum_{j=1}^{m}\frac{1}{T}\sum_{t=1}%
^{T}(z_{j,t}-\widehat{z}_{j,t}(\mathbf{f,\beta}_{j},\alpha_{j}\mathbf{))}%
^{2}=\sum_{j=1}^{m}\sum_{t=1}^{T}(z_{j,t}-\sum_{i=0}^{k}\beta_{j,i+1}%
f_{t+i}-\alpha_{j})^{2}. \label{ecmloss}%
\end{equation}
The optimal choices of $\mathbf{f=(}f_{1},...,f_{T+k})^{\prime}$ and
$\mathbf{\beta=(}\beta_{j,i})_{1\leq j\leq m,1\leq i\leq k+1}$
,\ $\mathbf{\alpha=(}\alpha_{1},...\alpha_{m})$ are given by
\begin{equation}
(\widehat{\mathbf{f}}\mathbf{,}\widehat{\mathbf{\beta}}\mathbf{)}=\arg
\min_{\mathbf{f,\beta,\alpha}}\text{MSE}(\mathbf{f,\beta,\alpha).} \label{SS}%
\end{equation}
Clearly if $\mathbf{f}$ is optimal, $\gamma\mathbf{f+}\delta$ is optimal too.
Thus, we can choose $\mathbf{f}$ so that $\sum_{t=1}^{T+k}f_{t}^{2}%
/(T+k)=1,$and $\sum_{t=1}^{T+k}f_{t}/(T+k)=0.$We call $\widehat{\mathbf{f}}%
\ $the first DPC of order $k$ of the observed series $\mathbf{z}%
_{1},...,\mathbf{z}_{t}.$ \ Note that the first DPC of order $0$ corresponds
to the first \ regular principal component of the data. Moreover $\ $the
matrix $\widehat{\mathbf{\beta}}$ \ \ contains the coefficients to be used to
reconstruct the $m$ series from $\widehat{\mathbf{f}}$ in an optimal way$.$\ 

Let $\mathbf{C}_{j}(\alpha_{j})=(c_{j,t,q}(\alpha_{j}))_{1\leq t\leq T+k,1\leq
q\leq k+1}$ \ be the $(T+k)\times(k+1)$ matrix defined by
\begin{equation}
c_{j,t,q}(\alpha_{j})=\left\{
\begin{array}
[c]{ccc}%
(z_{j,t-q+1}-\alpha_{j}) & \text{if } & 1\vee(t-T+1)\leq q\leq(k+1)\wedge t\\
0 & \text{if} & \text{otherwise}%
\end{array}
\right.  . \label{matC}%
\end{equation}
\ where $a\vee b=\max(a,b)\ $and $a\wedge b=\min(a,b).$ Let $\mathbf{D}%
_{j}(\mathbf{f,\beta}_{j}\mathbf{)}=(d_{j,t,q}(\mathbf{f,\beta}_{j}%
\mathbf{)})$ be the $(T+k)\times(T+k)$ given by%
\[
d_{j,t,q}(\mathbf{f,\beta}_{j}\mathbf{)}=\left\{
\begin{array}
[c]{ccc}%
\sum_{v=(t-k)\vee1}^{t\wedge T}\beta_{j,q-v+1}\beta_{j,t-v+1} & \text{if} &
(t-k)\vee1\leq q\leq(t+k)\wedge(T+k)\\
0 & \text{if} & \text{otherwise}%
\end{array}
\right.
\]
and
\begin{equation}
\mathbf{D}(\mathbf{f,\beta})=\sum_{j=1}^{m}\mathbf{D}_{j}(\mathbf{f,\beta}%
_{j}\mathbf{).} \label{matD}%
\end{equation}
Differentiating (\ref{ecmloss}) with respect to $f_{t}$ \ in Subsection
\ref{derg1} we get the following equation
\begin{equation}
\mathbf{f=D}(\mathbf{f,\beta})^{-1}\sum_{j=1}^{m}\mathbf{C}_{j}(\mathbf{\alpha
})\mathbf{\beta}_{j}\mathbf{.} \label{g1}%
\end{equation}
Obviously, the $\ $coefficients $\mathbf{\beta}_{j}$ and $\alpha_{j,}$ $1\leq
j\leq m,\ $\ can be obtained using the least squares estimator, that is
\begin{equation}
\left(
\begin{array}
[c]{c}%
\mathbf{\beta}_{j}\\
\alpha_{j}%
\end{array}
\right)  =\left(  \mathbf{F(f)}^{\prime}\mathbf{F(f)}\right)  ^{-1}%
\mathbf{F(f)}^{\prime}\mathbf{\ z}^{(j)}, \label{g2}%
\end{equation}
where $\mathbf{z}^{(j)}=(z_{j,1},...,z_{j,T})^{\prime}$ and $\mathbf{F(f)}$ is
the $T\times(k+2)$ matrix with $t$-th row ($f_{t},f_{t+1},...,f_{t+k},1).$
\ Then the first DPC is determined by equations (\ref{g1}) and (\ref{g2}%
).$\ \ $The second DPC \ is defined as the first DPC of the \ residuals
$r_{j,t}(\mathbf{f,\beta).}$ Higher order DPC are defined in a similar manner.
We will call $p$ the selected number of components.

To define an iterative algorithm to compute $(\widehat{\mathbf{f}}%
\mathbf{,}\widehat{\mathbf{\beta}}\mathbf{,}\widehat{\alpha})$ is enough to
\ give $\mathbf{f}^{(0)}$ and to describe how \ to compute $\mathbf{\beta
}^{(h)},\mathbf{\alpha}^{(h)},$ $f^{(h+1)}$ once $\mathbf{f}^{(h)}$ is known.
\ According to \ (\ref{g1}) and (\ref{g2}) \ \ a natural such a rule is given
by the following \ two steps:

\begin{description}
\item[step 1] Based on (\ref{g2}), define $\mathbf{\beta}_{j}^{(h)}$ and
$\alpha_{j}^{(h)},$ for $1\leq j\leq m$\ , by
\[
\left(
\begin{array}
[c]{c}%
\mathbf{\beta}_{j}^{(h)}\\
\alpha_{j}^{(h)}%
\end{array}
\right)  =\left(  \mathbf{F(f}^{(h)}\mathbf{)}^{\prime}\ \mathbf{F(f}%
^{(h)}\mathbf{)}\right)  ^{-1}\mathbf{F(f}^{(h)}\mathbf{)}^{\prime}%
\mathbf{z}^{(j)}.
\]

\item[step 2] Based on (\ref{g1}), define $\mathbf{f}^{(h+1)}$ by
\[
\mathbf{f}^{\ast}=\mathbf{D}(\mathbf{f}^{(h)}\mathbf{,\beta}^{(h)}%
\mathbf{,\alpha}^{(h)})^{-1}C(\mathbf{f}^{(h)}\mathbf{,\beta}^{(h)}%
\mathbf{,\alpha})\mathbf{\beta}^{(h)}%
\]
\ and \ \
\[
\mathbf{f}^{(h+1)}\mathbf{=(T+}k\mathbf{)}^{1/2}(\mathbf{f}^{\ast}%
-\overline{\mathbf{f}}^{\ast})\mathbf{/|||\mathbf{f}^{\ast}-\overline
{\mathbf{f}}^{\ast}||.}\bigskip
\]

\end{description}

The initial value $\mathbf{f}^{(0)}$ can be chosen equal to the standard (non
dynamic) first principal component, completed with $\ k$ zeros. The iterative
procedure is stopped when
\[
\frac{\text{MSE}(\mathbf{f}^{(h)}\mathbf{,\beta}^{(h)}\mathbf{,\alpha}%
^{(h)}\mathbf{)-}\text{MSE}(\mathbf{f}^{(h+1)}\mathbf{,\beta}^{(h+1)}%
\mathbf{,\alpha}^{(h+1)}\mathbf{)}}{\text{MSE}(\mathbf{f}^{(h)}\mathbf{,\beta
}^{(h)}\mathbf{,\alpha}^{(h)}\mathbf{)}}<\varepsilon
\]
for some value $\varepsilon.$

Note that we start with $m$ series of size $T.$ Assuming that we consider $p$
dynamic principal components let $\beta_{j,i,s}$ $1\leq j\leq m,1\leq i\leq
k+1,$ the coefficients $\beta_{j.i}$ corresponding to the $s-$th component,
$1\leq s\leq p.$ Then, the number of values required to reconstruct the
original series are the $(T+k)p$ values of the $p$ factors plus $(k+1)mp$
values for the coefficients $\beta_{j,i,s}$ plus the $m$ intercepts
$\alpha_{j}.$ Thus the proportion of the original information required to
reconstruct the series is $((T+k)p+(k+1)mp+m)/mT$ and when $T$ is large
compared to $k$ and $m$ is close to $p/m.$ In applications the number of lags
to reconstruct the series, $k,$ and the number of principal components, $p,$
need to be chosen. Of course the accuracy of the reconstruction improves when
any of these two numbers is enlarged, but also the size of the information
required will also increase. For large $T$ increasing the number of components
introduces more values to store than increasing the number of lags. However,
we should also take into account the reduction in MSE due to enlarging each of
these components. Is clear that increasing the number of lags after some point
will have a negligible effect on the reduction in MSE. Then, if the level of
the MSE is larger than desired, adding an additional component is call for.
Thus one possible strategy will be start with one factor and increase the
number of lags until the reduction of further lags is smaller than $\epsilon.$
Then a new factor is introduced and the same procedure is applied. The process
stops when the MSE\ reaches some satisfactory value. Note that this rule is
similar to what is generally used for determining the number $p$ in ordinary
principal components.

\section{Dynamic Principal Components when $k=1$\label{RECOn1}}

To illustrate the computation of the \ first DPC, let us consider the simplest
case $\ $of $k=1.$ Then, we \ search for $\widehat{\mathbf{\beta}}%
\mathbf{=}(\widehat{\beta}_{ji})_{1\leq j\leq m,1\leq i\leq2}$ and
$\widehat{\mathbf{f}}\mathbf{=(}\widehat{f}_{1},...,\widehat{f}_{T+1}%
)^{\prime}$ such that%
\begin{equation}
(\widehat{\mathbf{f}}\mathbf{,}\widehat{\mathbf{\beta}}\mathbf{)}=\arg\min
_{1}\sum_{t=1}^{T}\sum_{j=1}^{m}(z_{j,t}-\beta_{j,1}f_{t}-\beta_{j,2}%
f_{t+1})^{2}. \label{summ}%
\end{equation}
Put $a_{1}=\sum_{j=1}^{m}\beta_{j,1}^{2},$ $a_{2}=\sum_{j=1}^{m}\beta
_{j,2}^{2}$ and $b=\sum_{j=1}^{m}\beta_{j,1}\beta_{j,2},$then the matrix
$\mathbf{D=}\sum_{j=1}^{m}D_{j}$ defined in (\ref{matD}) can be written as
\[
D\mathbf{=}a_{2}\left(
\begin{array}
[c]{cccccc}%
a_{1}/a_{2} & b/a_{2} & 0 & 0 & ... & ...\\
b/a_{2} & 1+a_{1}/a_{2} & b/a_{2} & 0 & ... & ...\\
0 & b/a_{2} & 1+a_{1}/a_{2} & b/a_{2} & ... & ...\\
... & ... & ... & ... & ... & ...\\
0 & ... & ... & b/a_{2} & 1+a_{1}/a_{2} & b/a_{2}\\
0 & ... & ... & 0 & b/a_{2} & 1
\end{array}
\right)  .
\]
\ Let $\widehat{\mathbf{\beta}}^{(i)}=(\widehat{\beta}_{i,1},...,\widehat
{\beta}_{i,m}),i=$ $1,$ $2.$ It is shown in the appendix that \ if
$\widehat{\mathbf{\beta}}^{(1)}\neq\lambda\widehat{\mathbf{\beta}}^{(2)}$there
exists $|c|<1,\alpha,w_{1}$ and $w_{2}$ so that \
\begin{equation}
D\mathbf{=}\alpha\left(
\begin{array}
[c]{cccccc}%
w_{1} & -c & 0 & 0 & ... & ...\\
-c & 1+c^{2} & -c & 0 & ... & ...\\
0 & -c & 1+c^{2} & -c & ... & ...\\
... & ... & ... & ... & ... & ...\\
0 & ... & ... & -c & 1+c^{2} & -c\\
0 & ... & ... & 0 & -c & w_{2}%
\end{array}
\right)  . \label{matrixC}%
\end{equation}
\ Note that $\widehat{\mathbf{\beta}}^{(1)}=\lambda\widehat{\mathbf{\beta}%
}^{(2)}$ implies that putting $\widehat{\mathbf{f}}^{\ast}\mathbf{=(}%
\widehat{f}_{t}^{\ast})_{1\leq t\leq T}$ where $\widehat{f}_{t}^{\ast
}=\widehat{f}_{t}+\lambda\widehat{f}_{t+1}$we have%
\[
\sum_{t=1}^{T}\sum_{j=1}^{m}(z_{j,t}-\widehat{\beta}_{j,1}\widehat{f}%
_{t}-\widehat{\beta}_{j,2}\widehat{f}_{t+1})^{2}=\sum_{t=1}^{T}\sum_{j=1}%
^{m}(z_{j,t}-\widehat{\beta}_{j}\widehat{f}_{t}^{\ast}\ )^{2},
\]
and therefore, in this case the first DPC is as good for reconstructing the
series as the first classical PC.

Let $A_{0}$ be defined by%
\[
A_{0}\mathbf{=}\left(
\begin{array}
[c]{cccccc}%
1 & -c & 0 & 0 & ... & ...\\
-c & 1+c^{2} & -c & 0 & ... & ...\\
0 & -c & 1+c^{2} & -c & ... & ...\\
... & ... & ... & ... & ... & ...\\
0 & ... & ... & -c & 1+c^{2} & -c\\
0 & ... & ... & 0 & -c & 1
\end{array}
\right)  ,
\]
put $m_{1}=w_{1}-1,m_{2}=w_{2}-1$ and $\ $\ let $G=(G_{1},G_{2})$ be the
$(T+1)\times2$ dimensional matrix where $G_{1}^{\prime}=(m_{1}^{1/2},0,...,0)$
and $G_{2}^{\prime}=(0,...,0,m_{1}^{1/2}).$ We can write $D=\alpha
(A_{0}+GG^{\prime}),$and then according to the Proposal A.3.3 of Seber (1984)
we have%
\begin{align}
D^{-1}  &  =\frac{1}{\alpha}\left(  A_{0}^{-1}-A_{0}^{-1}G(I+G^{/}A_{0}%
^{-1}G)^{-1}G^{\prime}A_{0}^{-1}\right) \label{dd-1}\\
&  =\frac{1}{\alpha}(A_{0}^{-1}-A_{0}^{-1}GHG^{\prime}A_{0}^{-1}),\nonumber
\end{align}
where $H=(I+G^{/}A_{0}^{-1}G)^{-1}=(h_{i,h})$ is a $2\times2$ matrix. \ We
also have that \ $A_{0}^{-1}$ is of the form%
\begin{equation}
(A_{0}^{-1})_{i,h}=\frac{1}{1-c^{2}}c^{|i-h|}. \label{A0-1}%
\end{equation}
and then we get%
\[
A_{0}^{-1}GH=\frac{1}{1-c^{2}}\left(
\begin{array}
[c]{cc}%
m_{1}^{1/2}h_{11}+m_{2}^{1/2}c^{T}h_{21} & m_{1}^{1/2}h+m_{2}^{1/2}c^{T}%
h_{22}\\
& \ \\
m_{1}^{1/2}h_{11}c^{i-1}+m_{2}h_{21}c^{T-i+1} & m_{1}^{1/2}h_{12}c^{i-1}%
+m_{2}h_{22}c^{T-i+1}\\
& \\
\  & \ \\
m_{1}^{1/2}h_{11}c^{T}+m_{2}h_{21} & m_{1}^{1/2}h_{12}c^{T}+m_{2}h_{22}%
\end{array}
\right)
\]
and%
\begin{align}
(A_{0}^{-1}GHG^{\prime}A_{0}^{-1})_{ih}  &  =1/(1-c^{2})^{2}\ [(m_{1}%
^{1/2}h_{11}^{\ }c^{i-1}+m_{2}^{1/2}h_{21}c^{T-i+1})m_{1}^{1/2}c^{h-1}%
\nonumber\\
&  +(m_{1}^{1/2}h_{12}c^{i-1}+m_{2}h_{22}c^{T-i+1})m_{2}^{1/2}c^{T-h+1}%
]\nonumber\\
&  =A_{1}c^{i+h-2}+A_{2}c^{T-i+h}+A_{3}c^{2T-i-h+2}. \label{eslle}%
\end{align}

By \ (\ref{g1}) \ we have $\widehat{\mathbf{f}}\mathbf{=D}^{-1}\sum_{j=1}%
^{m}\mathbf{C}_{j}\widehat{\mathbf{\beta}}_{j},$where $\widehat{\mathbf{\beta
}}_{j}$ is given by (\ref{g2}) and $C_{j}=(Z_{1},Z_{2})$ where $Z_{1}^{\prime
}$ $=(z_{j,1},...,z_{j,T},0)$ and $Z_{2}^{\prime}$ $=(0,z_{j,1},...,z_{j,T}).$
Therefore, by (\ref{dd-1}), (\ref{A0-1}) and (\ref{eslle}) we obtain%

\[
\widehat{f}_{t}=\frac{1}{\alpha}\left[
{\displaystyle\sum\limits_{j=1}^{m}}
\widehat{\beta}_{j,1}%
{\displaystyle\sum\limits_{q=1}^{T}}
c^{|t-q|}z_{j,q}+%
{\displaystyle\sum\limits_{j=1}^{m}}
\widehat{\beta}_{j,2}%
{\displaystyle\sum\limits_{q=2}^{T+1}}
c^{|t-q|}z_{j,q-1}\right]  +R_{t}%
\]
where $R_{t}\rightarrow0$ except for $t$ close to $1$ or to $T.$

Suppose now that $\mathbf{z}_{t}$ is stationary, then except in both ends
$\widehat{f}_{t}$ can be approximated by the stationary process%
\[
\widehat{f}_{t}^{\ast}=\frac{1}{\alpha}\left[
{\displaystyle\sum\limits_{j=1}^{m}}
\widehat{\beta}_{j,1}%
{\displaystyle\sum\limits_{q=-\infty}^{\infty}}
c^{|t-q|}z_{j,q}+%
{\displaystyle\sum\limits_{j=1}^{m}}
\widehat{\beta}_{j,2}%
{\displaystyle\sum\limits_{q=-\infty}^{\infty}}
c^{|t-q|}z_{j,q-1}\right]  ,
\]
and the DPC is approximated as linear combinations of the geometrically and
symmetrically filtered \ series $z_{j,t}+%
{\displaystyle\sum\limits_{i=1}^{\infty}}
c^{i}(z_{j,t+i}+z_{j,t-i}),$ and $z_{j,t-1}+%
{\displaystyle\sum\limits_{i=1}^{\infty}}
c^{i}(z_{j,t-1+i}+z_{j,t-1-i}),1\leq j\leq m.$ These series give the largest
weight to the periods $\ t$ and $\ t-1$ respectively and the weights decrease
geometrically when we move away of these values. We conjecture that in the
case of the first DPC of order $k,$ a similar \ approximation \ outside both
ends of $\widehat{f}_{t}$ by an stationary process can be obtained.

\section{Monte Carlo simulation and two real examples\label{simul}}

We perform a Monte Carlo study \ using \ as vector series $\ \mathbf{z}%
_{t}=(z_{1,t},z_{2,t},z_{3,t})^{\prime}$, $1\leq t\leq T$ \ generated as
follows: let $v_{t},1\leq t\leq T+2,$ $w_{i,t},$ $1\leq i\leq3,1\leq t\leq T$,
i.i.d random variables with distribution N$(0,1),$ then $z_{i,t}%
=v_{t+i-1}+0.1w_{i,t},$ $1\leq i\leq3,$ $1\leq t\leq T$. We compute three
different principal components: \ (i) The ordinary principal component
(OPC)$,$ (ii) the dynamic principal component (DPC$_{k})$ proposed \ here with
$\ k,1,5$ and $10,$ (iii) Brillinger dynamic principal \ components
(BDPC$_{M}$) adapted for finite samples as follows:
\begin{equation}
f_{t}=%
{\displaystyle\sum\limits_{k=(-M)\vee(t-T)}^{M\wedge(t-1)}}
\mathbf{c}_{k}^{\prime}\mathbf{z}_{t-k},
\end{equation}
where $\mathbf{c}_{k}$ are the coefficients defined below (\ref{BR2}) in
Section 1. \ The values of $M$ where taken 10, 20 and 50. To reconstruct the
original series with the OPC we used $\ k=1,5,10$ lags \ and the corresponding
coefficients were obtained \ using least squares. \ To reconstruct the series
with DPC$_{k}$ we proceed as described in Section \ref{recons}. Finally, \ the
original series $\mathbf{z}_{t}$ \ were reconstructed using the BDPC$_{M}$
\ by%
\[
\widehat{z}_{i,t}=\sum_{j=(-M)\vee(-t+1)}^{M\wedge(T-t)}\mathbf{\beta}%
_{i,,j}f_{t+j}%
\]
where the $\mathbf{\beta}_{i,,j}$ are described below (\ref{BR2}) in Section
1. The cross spectrum matrix was computed using the function \texttt{mvspec}
in the ASTSA package \ with the \ R software. \ We took two values of $T:$ 100
and 500 and we make 500 replications. Table \ref{TableKeySim} shows the MSE of
\ the prediction residuals obtained with OPC, DPC$_{k}$ and BDPC$_{M}.$ We
observe that the procedure DPC$_{k}$ proposed here produces a much better
reconstruction of the original series than the OPC$_{k}$ and the BDPC$_{M}.$

\begin{center}
\bigskip%
\begin{table}[h] \centering
\begin{tabular}
[c]{llllllllll}\hline
$T$ & \multicolumn{3}{c}{OPC$_{k}$} & \multicolumn{3}{c}{DPC$_{k}$} &
\multicolumn{3}{c}{BDPC$_{M}$}\\\cline{2-10}
& \multicolumn{3}{c}{$k$} & \multicolumn{3}{c}{$k$} & \multicolumn{3}{c}{M}%
\\\cline{2-10}
& 1 & 5 & 10 & 1 & 5 & 10 & 10 & 20 & 50\\
100 & 1.31 & 0.78 & 0.67 & 0.89 & 0.018 & 0.016 & 2.05 & 2.08 & 2.17\\
500 & 1.42 & 0.79 & 0.66 & 0.97 & 0.034 & 0.025 & 2.03 & 2.03 & 2.03\\\hline
\end{tabular}
\caption{Mean Square Errors obtained in the Monte Carlo study.}\label{TableKeySim}%
\end{table}%

\end{center}

\subsection{Example 1}

\begin{figure}
[ptb]
\begin{center}
\includegraphics[
height=4.715in,
width=5.1465in
]%
{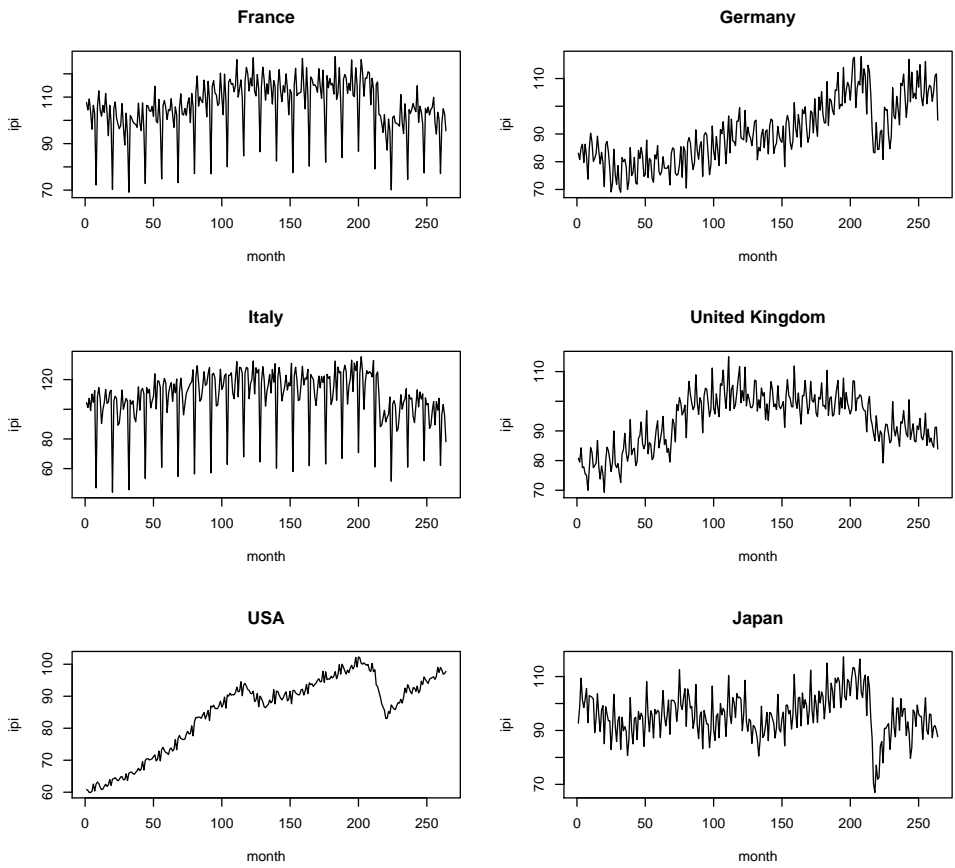}%
\caption{ Industrial production Index of six countries 1991-2012}%
\label{ipiserplot}%
\end{center}
\end{figure}

We use six series corresponding to \ the Industrial Production Index\ (IPI) of
France, Germany, Italy, United Kingdom, USA and Japan. We use monthly data
from January 1991 to December 2012 and the data are taken from Eurostat. The
seven series are plotted in Figure \ref{ipiserplot}.%

Let $\mathbf{f}_{k},$ $k\geq0,$ the first DPC$_{k}$. \ In Table \ref{Ipi72} we
show the percentage of\ variability explained by $\ \mathbf{f}_{0}$ and
$\mathbf{f}_{k}$ using $k$\ lags, computed as $EV_{j,k}=\min_{\mathbf{\beta
,\alpha}}$MSE$_{k}(\mathbf{f}_{j}\mathbf{,\beta,\alpha)/}\sum_{i=1}^{6}V_{i},$
for $j=0,k$ where $V_{i}$ is the variance of the series$\ i.$

\begin{center}%
\begin{table}[h] \centering
\begin{tabular}
[c]{lll}\hline
$k$ & $EV_{0,k}$ & $EV_{k,k}$\\\hline
$0$ & $63.07$ & $63.07$\\
$1$ & $66.19$ & $82.47$\\
$5$ & $76.66$ & $90.05$\\
$10$ & $77.98$ & $94.81$\\
$12$ & $80.00$ & $96.67$\\\hline
\end{tabular}
\caption{Explained variability  of the IPI series using the OPC and DPC   with  different number of lags }
\label{Ipi72}%
\end{table}%

\end{center}

We note that the reconstruction of the series \ using the DPC\ \ is notably
better that the one \ obtained by means of the OPC \ with the same lags.
\ Increasing the number of lags obviously improves the reconstruction obtained
by both components, although the improvement is larger with the DPC. With 12
lags the reconstruction error with the first DPC\ is smaller then 3.5\%. Table
\ \ref{table coef} includes the coefficients of the six \ IPI series in the
ordinary PC and in the first DPC with $k=1.$

\begin{center}%
\begin{table}[h] \centering
\begin{tabular}
[c]{lllll}\hline
PC & PC(0) & PC(1) & DPC(0) & DPC(1)\\\hline
-0.456 & -0.456 & -0.001 & -3.951 & 3.965\\
-0.285 & -0.275 & -0.034 & -1.509 & 1.492\\
-0.719 & -0.750 & 0.099 & -6.548 & 6.577\\
-0.298 & -0.269 & -0.092 & -2.114 & 2.111\\
-0.241 & -0.198 & -0.138 & -0.787 & 0.760\\
-0.212 & -0.212 & -0.001 & -1.885 & 1.894\\\hline
\end{tabular}
\caption{Coefficients to reconstruct the IPI series by uing OPC and DPC with one lag }
\label{table coef}%
\end{table}%
\smallskip
\end{center}

\begin{figure}
[ptb]
\begin{center}
\includegraphics[
height=5.5971in,
width=6.109in
]%
{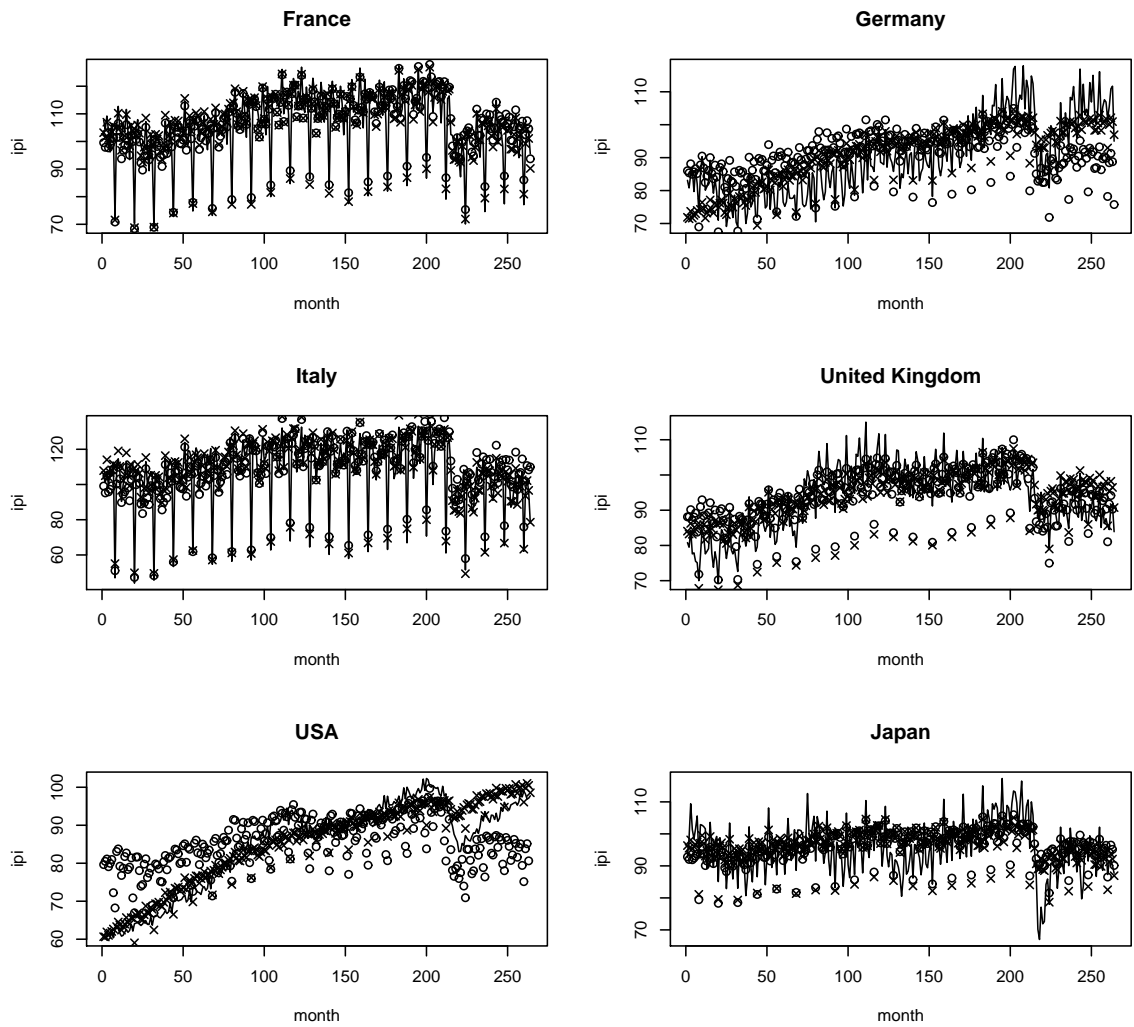}%
\caption{ Values of the original and reconstructed series of Example 1 with
the OPC (o) and DPC (*) with one lag}%
\label{predipiplot}%
\end{center}
\end{figure}

For the OPC the coefficients in the first column in Table \ref{table coef}
coincide with the weights given to each country in the definition of the OPC.
Thus, the first OPC\ gives the largest weight to Italy and then France,
because of the strong seasonality of these series which have the largest
variability. The second and third columns show that for reconstructing the
original variables including the lag of the OPC is practically irrelevant. The
fourth and fifth columns show that the DPC with one lag is almost equivalent
to using the first difference of the DPC in the reconstruction of the series.

\ Figure \ref{predipiplot} shows the original and reconstructed values using
the \ first OPC and the first DPC, both with one lag. We can see that the
reconstruction obtained with the DPC is clearly better than the one obtained
with the OPC for Germany and USA. In the other cases the reconstruction with
the DPC is still better but the differences are smaller and therefore more
difficult to detect in the plots.%

\begin{figure}
[ptb]
\begin{center}
\includegraphics[
height=5.7086in,
width=6.231in
]%
{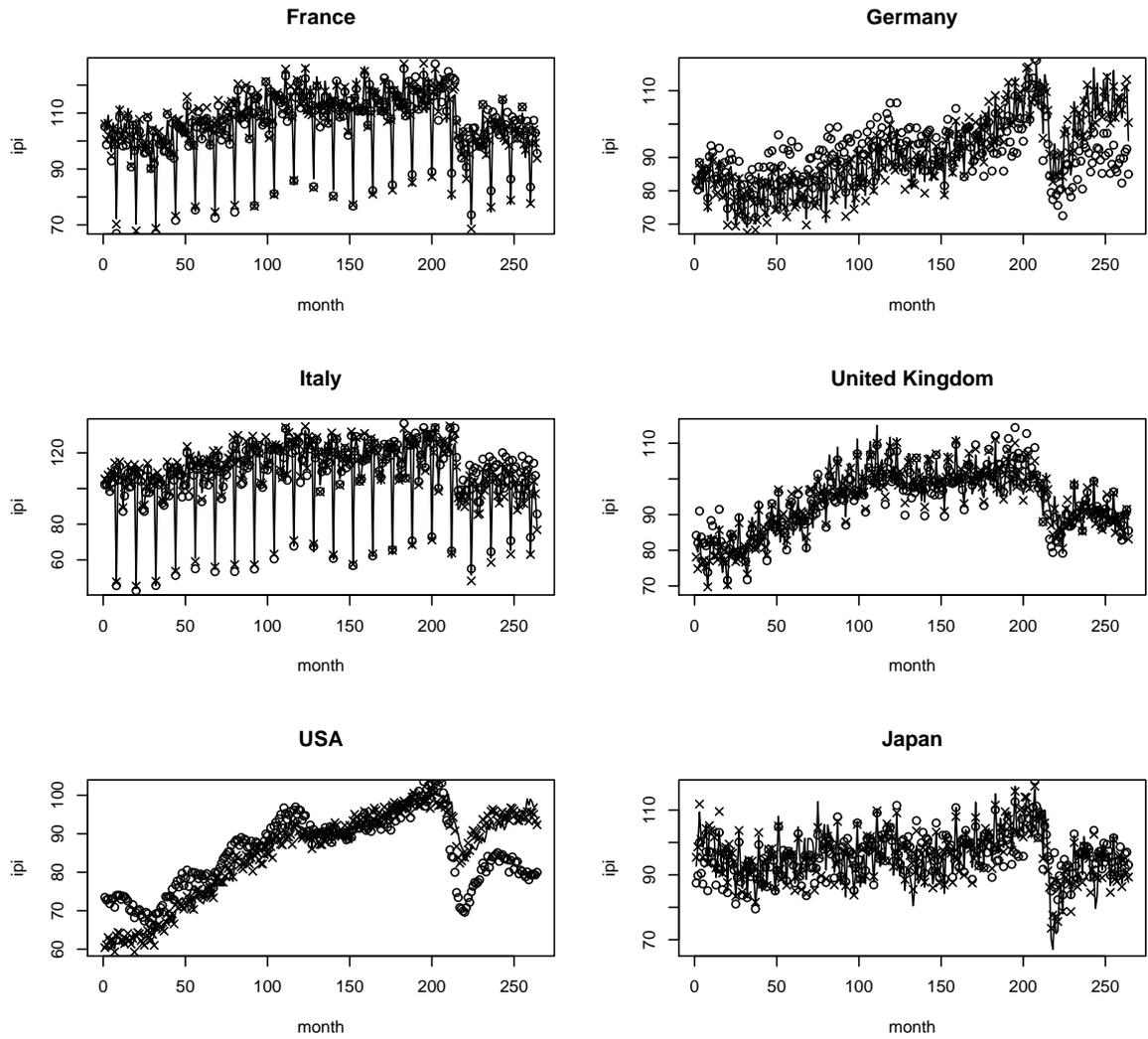}%
\caption{Values of the original and reconstructed series with the OPC (o) and
DPC (*) with twelve lags}%
\label{pred12fig}%
\end{center}
\end{figure}

Figure \ref{pred12fig} is similar to figure\ \ref{predipiplot}\ but with
twelve lags. \ Note that the reconstruction errors are significantly smaller
\ than in the case of one lag, and that there is an important improvement \ of
the reconstruction series when using the DPC \ instead of the OPC.%

\subsection{Example 2.}

In this example the data set is composed of 31 \ daily stock prices in the
stock market in Madrid corresponding to the 251 trading days of the year 2004.
These 31 series are the main components of the IBEX (general index of the
Madrid stock market). The source of the data is the Ministry of Economy,
Spain. In Table \ref{TableIBEX} we show the explained variability of the
reconstructed series \ using the DPC and OPC with different lags

\begin{center}%
\begin{table}[h] \centering
\begin{tabular}
[c]{lll}\hline
$k$ & $EV_{0,k}$ & $EV_{k,k}$\\\hline
0 & 0.598 & 0.598\\
1 & 0.602 & 0.822\\
5 & 0.610 & 0.873\\
10 & 0.620 & 0.881\\\hline
\end{tabular}
\caption{ Explained variability  of the OPC and DPC  for the stock prices
series with different number of lags}\label{TableIBEX}%
\end{table}%
\bigskip
\end{center}

%

\begin{figure}
[ptb]
\begin{center}
\includegraphics[
height=5.58in,
width=6.091in
]%
{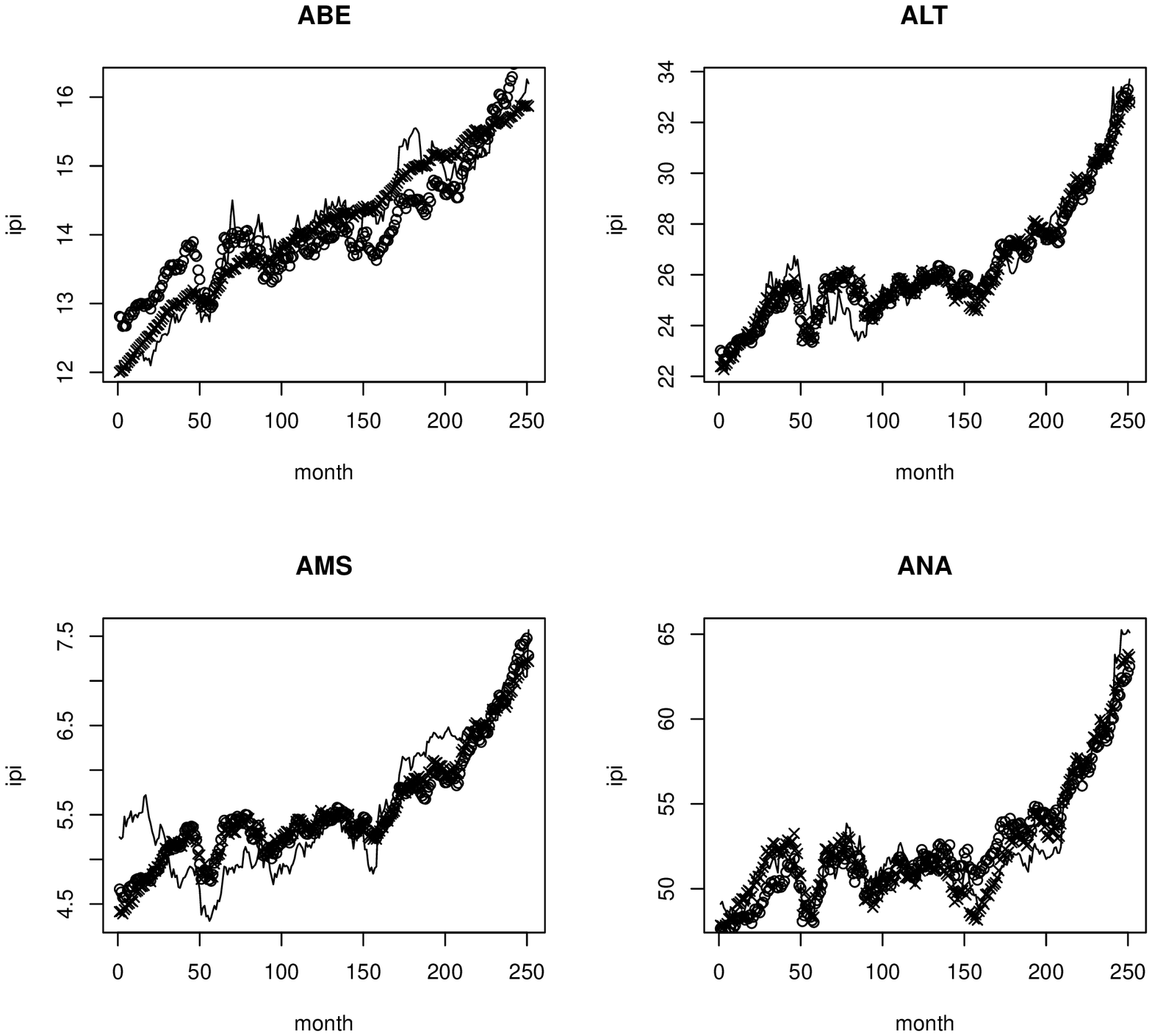}%
\caption{ Values of the original and reconstructed of the first four stocks
chosen in alphabetic orders. The reconstruction was made with the OPC (o) and
DPC (*) using one lag}%
\label{figIBEXpr1}%
\end{center}
\end{figure}
\begin{figure}
[ptb]
\begin{center}
\includegraphics[
height=4.5173in,
width=4.9306in
]%
{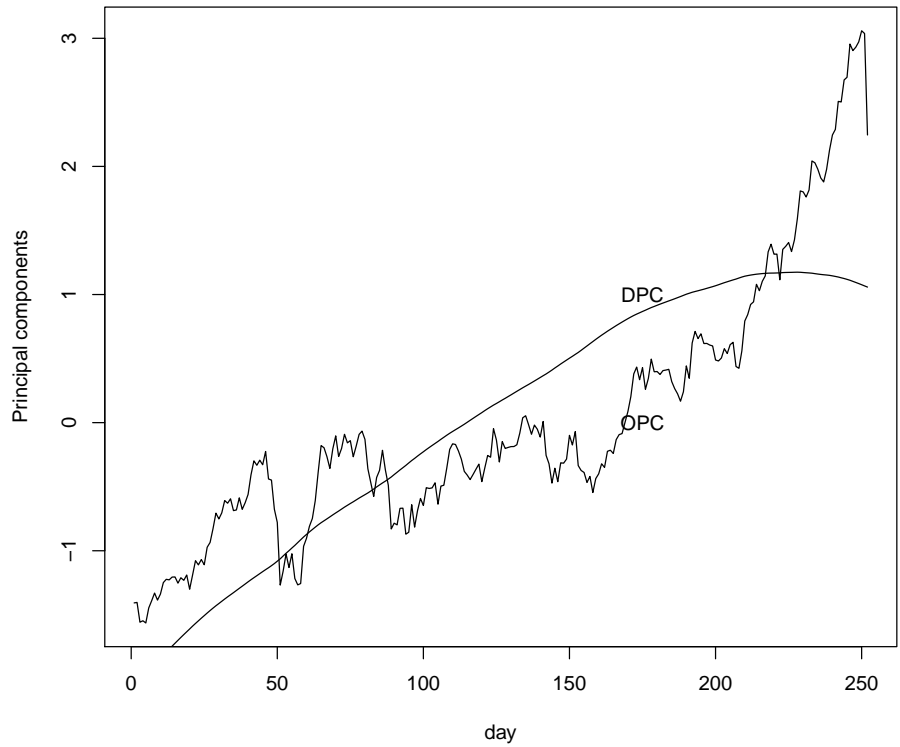}%
\caption{First OPC and DPC for the stock prices series }%
\label{figIBEXcomp}%
\end{center}
\end{figure}

In Figure \ref{figIBEXpr1} we show the first four series in alphabetic order
out of the thirty one \ and their reconstruction obtained by the first OPC and
DPC with one lag. As shown in Table \ref{TableIBEX} including one lag in the
OPC does not make much difference in the results, but it has a deep effect
when using the DPC. In fact, in the case of the DPC, the coefficient of the
one lag variable is very close but with opposite sign to the instantaneous
coefficient and therefore \ the reconstruction is \ similar to the one
obtained using the first difference of the first DPC without lags.\ Figure
\ref{figIBEXcomp} presents the first OPC and the DPC. The dynamic principal
components seems to be very useful to represent the general trend of the set
of time series.%

\section{Robust Dynamic Principal Components \label{SecRDC}}

As most of the procedures minimizing the mean square error, the DPC defined by
(\ref{SS}) is not robust. In fact a \ very small fraction of outliers may have
an unbounded influence on $(\mathbf{f,\alpha,\beta).}$ For this reason we are
going to study a robust alternative. One of the standard procedures to obtain
robust estimates for many statistical models is to replace the minimization
\ of the mean square scale for the minimization of a robust M-scale. This
strategy was used for many statistical models, including among other linear
regression (Rousseeuw \ and Yohai, 1984), the estimation of a scatter matrix
and multivariate location \ for multivariate data \ (Davis, 1987) and to
estimate the ordinary principal components (Maronna, 2005). The estimators
defined by means of a robust M-scale are called S-estimators. \ In this
section we extend the S-estimators for the case of the DPC.

Special care is required for time series with strong seasonality. The reason
is that a robust procedure may take the values corresponding to a particular
season \ which is \ very different to the others \ as outliers, \ and
therefore downweight these values. As a consequence, the \ reconstruction of
\ these observations may be affected by large errors. Thus, the procedure we
present here assumes that the series have been adjusted by seasonality and
therefore this problem is not present.

\subsection{S-Dynamic Principal Components}

Let $\rho_{0}$ be a symmetric, non-decreasing function for $x\geq0$ and
$\rho_{0}(0)=0.$Given a sample $\mathbf{x=(}x_{1},...,x_{n})$, the M-scale
estimator $S(\mathbf{x)}$ is defined as the value $s$ solution of%
\begin{equation}
\frac{1}{n}%
{\displaystyle\sum\limits_{i=1}^{n}}
\rho_{0}\left(  \frac{x_{i}}{s}\right)  =b. \label{M-sc}%
\end{equation}
\bigskip\ \ If $\rho_{0}$ is bounded, then the breakdown point to $\infty$\ of
$S(\mathbf{x),}$ that is, the minimum fraction \ of outliers \ than can take
$S(\mathbf{x)}$\ to $\infty$ is $b/\max\rho_{0}.$ Moreover, the breakdown
point to 0, that is, \ the minimum fraction of \ \ inliers that can take
$S(\mathbf{x)}$ to 0, \ is $1-$ $(b/\max\rho_{0}$). Note that if $b/\max
\rho_{0}=0.5$ both breakdown points are 0.5 (see section 3.2.2. in Maronna,
Martin and Yohai, 2006). \ In what follows we assume without loss of
generality that $\max\rho_{0}=1$. We also assume that $b=0.5$ so that both
breakdowns are equal 0.5. Moreover $\rho_{0}$ is chosen so that $E_{\phi}%
(\rho_{0}(\mathbf{x))=}b\mathbf{,}$where $\phi$ is the standard normal
distribution. \ This condition guarantees that for normal samples
$S(\mathbf{x)}$ is a consistent estimator of the standard deviation. One very
popular family of $\rho$ functions is the Tukey biweight family defined by%
\[
\rho_{c}^{T}(x)=\left\{
\begin{array}
[c]{ccc}%
1-\left(  1-(x/c)^{2}\right)  ^{3} & \text{if} & |x|\leq c\\
1 & \text{if} & |x|>c
\end{array}
\right.  .
\]

Then,\ we can define the first S-DPC as follows:\ for $\ 1\leq j\leq m,$ let
$\mathbf{r}_{j}\mathbf{(f,\beta}_{j},\alpha_{j}\mathbf{)=}$\newline%
$\mathbf{(}r_{j,t}(\mathbf{f,\beta}_{j},\alpha_{j}\mathbf{)})_{1\leq t\leq
T},$ where $r_{j,t}(\mathbf{f,\beta}_{j},\alpha_{j}\mathbf{)=}z_{j,t}%
-\sum_{i=0}^{k}\beta_{j,i}f_{t+i}-\alpha_{j}.$\ Define $\ $%
\begin{equation}
\text{SRS}(f,\mathbf{\beta,\alpha)=}\sum_{j=1}^{m}S^{2}(\mathbf{r}%
_{j}\mathbf{(f,\beta}_{j},\alpha_{j}\mathbf{)),} \label{SRS}%
\end{equation}%
\begin{equation}
(\widehat{\mathbf{f}},\widehat{\mathbf{\beta}},\widehat{\mathbf{\alpha}%
}\mathbf{)=}\arg\min_{\mathbf{f,\beta}}\text{SRS}(f,\mathbf{\beta,\alpha),}
\label{lossrob}%
\end{equation}
then \ $\widehat{\mathbf{f}}$ is the the first S-DPC and $\widehat
{\mathbf{\beta}}$ $\ $\ and $\widehat{\mathbf{\alpha}}$ \ are the coefficients
\ to reconstruct the $z_{j,t}$'s from $\widehat{\mathbf{f}}.$

Note that the only difference with the definition given in (\ref{SS}) is that
instead of minimizing the MSE of the residuals, we minimize \ the sum of
squares of the robust M-scales applied to the residuals of the $\ m$
series.\ \ Put $\psi=\rho^{\prime},w(u)=\psi(u)/u,$
\begin{equation}
s_{j}=s_{j}(\mathbf{f,\beta}_{j},\alpha_{j})=S(\mathbf{r(f,\beta}%
_{j}\mathbf{,}\alpha_{j}\mathbf{)).} \label{fineq1}%
\end{equation}
Note that $s_{j}$ satisfies%
\begin{equation}
\frac{1}{T}\sum_{t=1}^{T}\rho\left(  \frac{z_{v-}\sum_{i=0}^{k}\beta
_{j,i+1}f_{v+i}-\alpha_{j}}{s_{j}}\right)  =b. \label{sj}%
\end{equation}
Define the weights
\begin{equation}
w_{j,t}=w_{j,t}(\mathbf{f,\beta}_{j},\alpha_{j})=w_{0}\left(  \frac
{r_{j,t}(\mathbf{f,\beta}_{j})}{s_{j}}\right)  ,\text{ }1\leq j\leq m,\text{
}1\leq t\leq T \label{fineq2}%
\end{equation}
and%

\begin{equation}
W_{j,t,v}=W_{j,t,v}(\mathbf{f,\beta},\mathbf{\alpha},s)=\frac{s_{j}^{2}%
w_{j,v}(\mathbf{f,\beta}_{j},\alpha_{j},s_{j})}{\sum_{h=(t-k)\vee1}^{t\wedge
T}w_{j,h}(\mathbf{f,\beta}_{j},\alpha_{j},s_{j})r_{j,h}^{2}}\ , \label{WW}%
\end{equation}
where $\mathbf{s=}(s_{1},...s_{m})$. Let $\mathbf{C}_{j}(\mathbf{f,\beta}%
_{j},s)=(c_{j,t,q}(\mathbf{f,\beta}_{j},s))_{1\leq t\leq T+k,0\leq q\leq k}$
\ be the\newline$(T+k)\times(k+1)$ matrix defined by
\begin{equation}
c_{j,t,q}(\mathbf{f,\beta},\mathbf{\alpha,s})=\left\{
\begin{array}
[c]{ccc}%
W_{j,t-q+1}(\mathbf{f,\beta},\mathbf{\alpha},\mathbf{s})(z_{j,t-q+1}%
-\alpha_{j}) & \text{if } & 1\vee(t-T+1)\leq q\leq(k+1)\wedge t\\
0 & \text{if} & \text{otherwise}%
\end{array}
\right.  , \label{CC}%
\end{equation}
$\mathbf{D}_{j}$ $(\mathbf{f,\beta},\mathbf{\alpha,s)}=(d_{j,t,q}%
(\mathbf{f,\beta},\mathbf{\alpha,s)})$ the $(T+k)\times(T+k)$ matrix \ with
elements%
\[
d_{j,t,q}(\mathbf{f,\beta},\mathbf{\alpha,s)}=\left\{
\begin{array}
[c]{ccc}%
\sum_{v=(t-k)\vee1}^{t\wedge T}W_{j,t,v}\beta_{j,q-v+1}\beta_{j,t-v+1} &
\text{if} & (t-k)\vee1\leq q\leq(t+k)\wedge(T+k)\\
0 & \text{if} & \text{otherwise}%
\end{array}
\right.
\]
\ \ and
\begin{equation}
\mathbf{D}(\mathbf{f,\beta,\alpha,s})=\sum_{j=1}^{m}\mathbf{D}_{j}%
(\mathbf{f,\beta},\mathbf{\alpha,s).} \label{DD}%
\end{equation}
\ Differentiating (\ref{sj}) \ with respect to $f_{t}$ we get the following
equation
\begin{equation}
\mathbf{f=D}(\mathbf{f,\beta,\alpha},\mathbf{s})^{-1}\sum_{j=1}^{m}%
\mathbf{C}_{j}(\mathbf{f,\beta,\alpha,s})\mathbf{\beta}_{j}\mathbf{.}
\label{impa5}%
\end{equation}
Let $\mathbf{F(f)}$ be the $T\times(k+2)$ matrix with $t$-th row
($f_{t},f_{t+1},...,f_{t+k},1)$ and $W_{j}(\mathbf{f,\beta,}s)$ be the
diagonal matrix with diagonal equal to $w_{j,1}((\mathbf{f,\beta}%
_{j},s),...,w_{j,T}(\mathbf{f,\beta}_{j},s)$. Then differentiating
\ (\ref{sj}) with respect to $\beta_{j,i}$ \ and $\alpha_{j}$we get
\begin{equation}
\left(
\begin{array}
[c]{c}%
\mathbf{\beta}_{j}\\
\alpha_{j}%
\end{array}
\right)  =\left(  \mathbf{F(f)}^{\prime}W_{j}(\mathbf{f,\beta}_{j}%
,s)\mathbf{F(f)}\right)  ^{-1}\mathbf{F(f)W}_{j}(\mathbf{f,\beta}%
_{j},s)^{\prime}\mathbf{z}^{(j)}. \label{impa6}%
\end{equation}
Then the first S-PDC is determined by equation (\ref{fineq1}),(\ref{impa5})and
(\ref{impa6}).$\ $Note that \ the estimator defined by (\ref{SS}) is an
S-estimate corresponding to $\rho_{0}^{2}(u)=u^{2}$ and $b=1.$ In this case
$w(u)=2$ ad then we have $w_{j,v}=1$ and $W_{j,v}=T$ for all $j$ and all$v.$
Then for this case (\ref{impa5}) and (\ref{impa6}) become (\ref{g1}) and
(\ref{g2}) respectively.

The second S-DPC \ is defined as the first S-DPC of the \ residuals
$r_{j,t}(\mathbf{f,\beta).}$ Higher order S-DPC are defined in a similar manner.

One important point is the choice of $b.$ At first sight, $b=.5$ may seem a
good choice, since in this case we are protected against up to 50 \% of large
outliers. However, the following argument shows that this choice may not be
convenient. \ The reason is that \ with this choice, the procedure has the so
called 50\% exact fitting property. This means that when 50 \% of \ the
$r_{j,t}\mathbf{\mathbf{(f,\beta}_{j},\alpha_{j}\mathbf{)}}$\textbf{\textbf{s}%
} are zero the scale $S(\mathbf{r}_{j}\mathbf{(f,\beta}_{j},\alpha
_{j}\mathbf{))}$ is 0 no matter the value of the remaining values. \ Moreover,
if 50 \% of \ the $|r_{j,t}\mathbf{\mathbf{(f,\beta}_{j},\alpha_{j}%
\mathbf{)|}}$ are small the scale $S(\mathbf{r}_{j}\mathbf{(f,\beta}%
_{j},\alpha_{j}\mathbf{))}$ is small too. Then when $b=0.5,$ the procedure may
choose $\ \mathbf{f,\beta}$ and $\mathbf{\alpha}$ so to reconstruct the values
\ corresponding to 50\% of the periods even if the dataset do not contain
outliers.. For this reason it is convenient to choose a smaller value as $b,$
as for example $b=.10.$ In that case $\ $to obtain $S(\mathbf{r}%
_{j}\mathbf{(f,\beta}_{j},\alpha_{j}\mathbf{))}=0,$ it is required \ that 90\%
of the $r_{j,t}\mathbf{\mathbf{(f,\beta}_{j},\alpha_{j}\mathbf{)}}%
$\textbf{\textbf{s}} be 0.

One may wonder why for regression \ is \ common to use $b=0.5$ and the 50\%
exact fitting property does not bring the problems mentioned above. The reason
is that in this case, \ if there are no outliers, the regression hyperplane
fitting 50\% of the observations also fits the remaining 50\%. This does not
occur in the case of the dynamic principal components.

\subsection{Computational algorithms for the S-dynamic principal \ components}

The compute the first S-DPC we propose to use an iterative algorithm. \ We
start the computing algorithm in step 0, and \ denote by \ $\mathbf{f}%
^{(h)},\mathbf{\beta}^{(h)}$ $\mathbf{\alpha}^{(0)}$and $\mathbf{s}$ \ the
\ values computed in step $\ h.$

\ The initial value $\mathbf{f}^{(0)}$ can be chosen equal to a regular (non
dynamic) robust \ principal component, for example the one proposed in Maronna
(2005). \ Once $\mathbf{f}^{(0)}$ \ is computed we can use this value to
compute a matrix $F^{(0)}=F$ with $i$-th row $(f_{i}^{(0)},f_{i+1}%
^{(0)},...,f_{i+k}^{(0)},1)$. \ The $\ j$-th $\operatorname{row}$ of
$\mathbf{\beta}^{(0)}$ \ and $\alpha_{j}^{(0)}$can be obtained using a
regression S-estimate taking $\mathbf{z}^{(j)}$ as response and $F^{(0)}$ as
design matrix. \ Finally $s_{j}^{(0)}=S(\mathbf{r}_{j}(\mathbf{f}%
^{(0)},\mathbf{\beta}^{(0)}).$

Then to define the algorithm is enough to describe how \ to compute
$(\mathbf{f}^{(h+1)},$\newline$\mathbf{\beta}^{(h+1)},s^{(h+1)})$ once
$(\mathbf{f}^{(h)},\mathbf{\beta}^{(h)},s^{(h)})$ is known. \ This is done in
the following three steps:

\begin{description}
\item[step 1] According to (\ref{impa5}),compute
\[
\mathbf{f}^{\ast}=\mathbf{D}(\mathbf{f}^{(h)}\mathbf{,\beta}^{(h)}%
\mathbf{,\alpha^{(h)},s}^{(h)})^{-1}C(\mathbf{f}^{(h)}\mathbf{,\beta}%
^{(h)}\mathbf{,\alpha}^{(h)}\mathbf{s}^{(h)})\mathbf{\beta}^{(h)}%
\]
\ \ \ and \ put \ $\mathbf{f}^{(h+1)}\mathbf{=(T+}k\mathbf{)}^{1/2}%
(\mathbf{f}^{\ast}-\overline{\mathbf{f}}^{\ast})\mathbf{/|||\mathbf{f}^{\ast
}-\overline{\mathbf{f}}^{\ast}||.}$

\item[step 2] By (\ref{impa5}), calling $\mathbf{W}_{j}^{(h)}=\mathbf{W}%
_{j}(\mathbf{f}^{(h)}\mathbf{,\beta}^{(h)},\mathbf{\alpha}^{(h)},s^{(h)})$
compute $\ $\ the $j$-th row by%
\[
\left(
\begin{array}
[c]{c}%
\mathbf{\beta}_{j}^{(h+1)}\\
\alpha_{j}^{(h+1)}%
\end{array}
\right)  =\left(  \mathbf{F(f}^{(h+1)}\mathbf{)}^{\prime}\mathbf{W}_{j}%
^{(h)}\mathbf{F(f}^{(h+1)}\mathbf{)}\right)  ^{-1}\mathbf{F(f}^{(h+1)}%
\mathbf{)W}_{j}^{(h)\prime}\mathbf{z}^{(j)}%
\]
for $1\leq j\leq m$.

\item[step 3] Compute $s_{j}^{(h+1)}=S(\mathbf{r}_{j}(\mathbf{f}%
^{(h+1)},\mathbf{\beta},\alpha_{h+1})).$
\end{description}

The procedure is stopped when%
\[
\frac{\text{SRS}(f^{(h)},\mathbf{\beta}^{(h)}\mathbf{,\alpha}^{(h)}%
\mathbf{)-}\text{SRS}(f^{(h+1)},\mathbf{\beta}^{(h+1)}\mathbf{,\alpha}%
^{(h+1}\mathbf{)}}{\text{SRS}(f^{(h)},\mathbf{\beta}^{(h)}\mathbf{,\alpha
}^{(h)}\mathbf{)}}<\varepsilon,
\]
where $\varepsilon$ is a fixed small value.

A procedure similar to the one described at the end of Section \ref{recons}
can be used to determine a convenient number of lags and components replacing
the MSE\ by the SRS.

\subsection{Example 3 \label{SectionRobEx}}

We will use the data of example 2 to illustrate the performance of the robust
DPC. This dataset was modified as follows: each of the 7781 values composing
the dataset was modified with 5\% probability adding 20 to the true value. In
Table \ref{TableIBEX CONT} we include MSE \ in the reconstruction of the
series with the DPC. Since the DPC is very sensitive to the presence of
outliers, we also compute the S-DPC. Since the MSE is \ very sensitive to
outliers, we evaluate \ the performance of the principal components to
reconstruct the series by using the SRS criterion. We take as $\rho$ the
bisquare function with $c=5.13$ and $\ b=0.1.$ These values make the M-scale
consistent to the standard deviation \ in the Gaussian case. \ Table
\ref{TableIBEX CONT} gives the MSE of the non DPC$_{k}$ and the SRS for the
DPC$_{k}$ and S-DPC$_{k}$ for $k=1,5$ and $10$.

\begin{center}%
\begin{table}[h] \centering
\begin{tabular}
[c]{llll}\hline
$k$ & MSE of the DPC$_{k}$ & SRS of the DPC$_{k}$ & SRS of the S-DPC$_{k}%
$\\\hline
1 & \multicolumn{1}{c}{309.70} & \multicolumn{1}{c}{106.69} &
\multicolumn{1}{c}{39.84}\\
5 & \multicolumn{1}{c}{295.84} & \multicolumn{1}{c}{119.03} &
\multicolumn{1}{c}{37.81}\\
10 & \multicolumn{1}{c}{274.74} & \multicolumn{1}{c}{111.33} &
\multicolumn{1}{c}{31.95}\\\hline
\end{tabular}
\caption{ MSE and SRS of the DPC$_k$ and S DPC$_k$  for the  contaminated stock prices
series}\label{TableIBEX CONT}%
\end{table}%

\end{center}

\begin{figure}
[ptb]
\begin{center}
\includegraphics[
height=5.5426in,
width=6.0502in
]%
{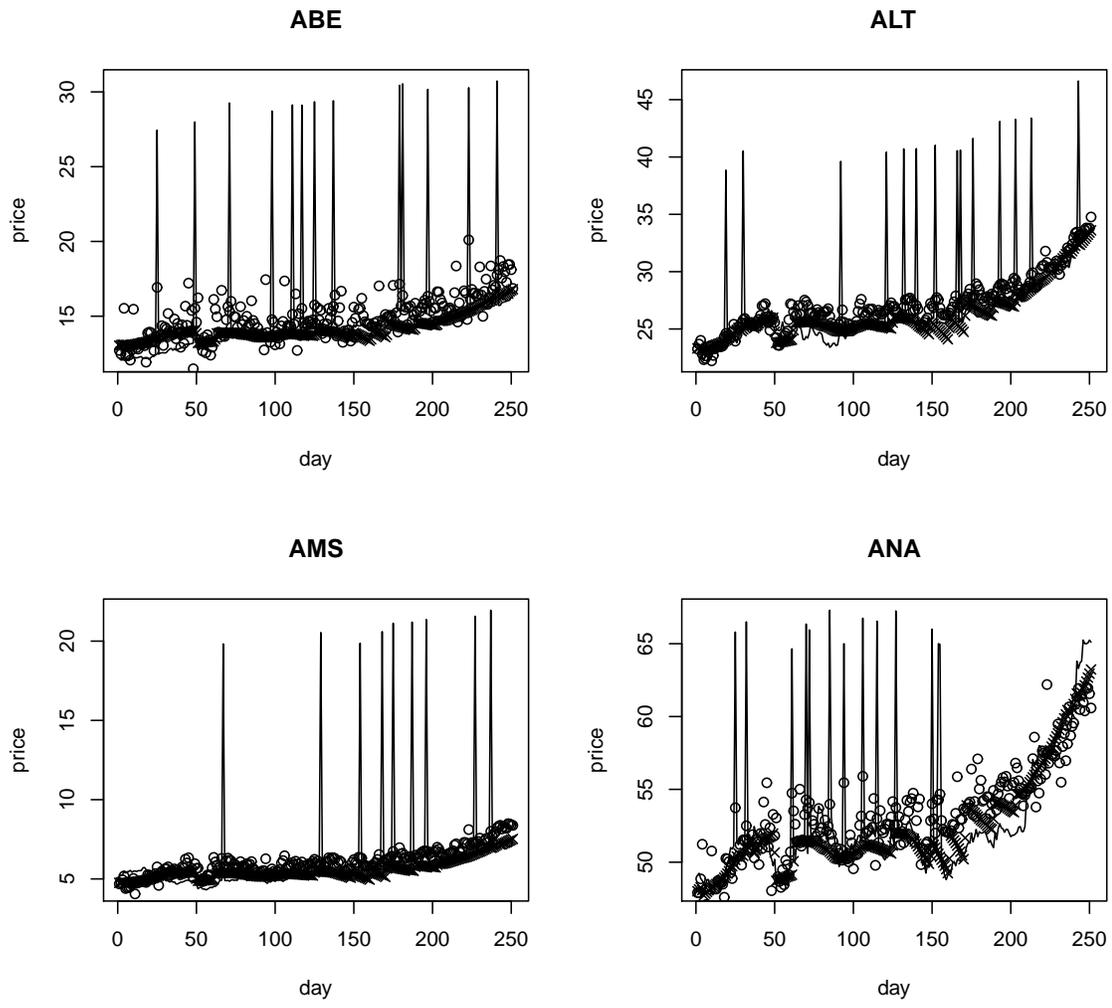}%
\caption{Contaminated Stock prices series and their reconstruction by DPC (o)
and by \ S-DPC (x)}%
\label{figure-pred-ibex-rob}%
\end{center}
\end{figure}

Figure \ref{figure-pred-ibex-rob} shows the reconstruction of the four stock
prices by using the DPC and the S-DPC. It can be seen, as expected, that the
robust methods has a better performance.%

\section{Conclusions \label{section: con}}

We have proposed two dynamic principal components procedures for multivariate
time series: the first one using a minimum squared error criterion to evaluate
the reconstruction of the original time series and the second one based on a
robust scale. These procedures, in contrast to previous ones, can also be
applied for nonstationary time series. A Monte Carlo study shows that the
proposed dynamic principal component based on the MSE criterion can improve
considerably the reconstruction obtained by both ordinary principal components
and a finite sample version of Brillinger approach. We have also shown in an
example that the robust procedure based on a robust scale is not much affected
by the presence of outliers.

A simple heuristic rule to determine a convenient value for the number of
components, $p,$ and the number of lags, $k,$ is suggested. However, further
research may lead to better methods to choose these parameters in order to
balance accuracy in the series reconstruction and economy in the number of
values stored for that purpose.

\section{Appendix \label{Sec: Ap}}

\subsection{Proof of (\ref{impa5})\label{derg1}}

Differentiating MSE$(\mathbf{f,\beta,\alpha)}$ with respect to $f_{t}$ for
$t=1,...,T+k\ $\ we get \newline$\sum_{j=1}^{m}\sum_{v=(t-k)\vee1}^{t\wedge
T}\left(  z_{j,v}-\sum_{i=0}^{k}\beta_{j,i+1}f_{v+i}\right)  \beta
_{j},_{t-v+1}=0,$where $a\wedge b$ \ denote minimum of $\ a$ and $b$ and
$a\vee b$ maximum. \ Then, we have
\begin{equation}
\sum_{j=1}^{m}\sum_{v=(t-k)\vee1}^{t\wedge T}(z_{j,v}-\alpha_{j}%
)\beta_{j,t-v+1}=\sum_{j=1}^{m}\sum_{v=(t-k)\vee1}^{t\wedge T}\sum_{i=0}%
^{k}\beta_{j,i+1}\beta_{j,t-v+1}f_{v+i} \label{1a}%
\end{equation}
that can be written as%
\begin{equation}
a_{t}(\mathbf{\beta)}=b_{t}(\mathbf{f,\beta),} \label{eqfun}%
\end{equation}
where $a_{t}(\mathbf{\beta)}$ and $b_{t}(\mathbf{f,\beta)}$ are the left and
right side of (\ref{1a}) respectively. Putting $q=t-v+1$ we have%
\begin{align}
a_{t}(\mathbf{\beta)}  &  =\sum_{j=1}^{m}\sum_{v=(t-k)\vee1}^{t\wedge
T}(z_{j,v}-\alpha_{j})\beta_{j,t-v+1}\label{lesi}\\
&  =\sum_{j=1}^{m}\sum_{q=1\vee(t-T+1)}^{(k+1)\wedge t}(z_{j,t-q+1}-\alpha
_{j})\beta_{j,q}.
\end{align}
and calling $\mathbf{a}(\mathbf{\beta})=(a_{1}(\mathbf{\beta}),...,a_{T+k}%
(\mathbf{\beta}))^{\prime}$
\begin{equation}
\mathbf{a}(\mathbf{\beta})=\sum_{j=1}^{m}\mathbf{C}_{j}(\alpha_{j}%
\mathbf{)\beta}_{j}\mathbf{.} \label{abeta}%
\end{equation}
where $C_{j}$ is given by (\ref{matC})

Now we will get an expression for $b_{t}(\mathbf{f,\beta).}$ Putting $q=v+i$
we get%

\[
b_{t}(\mathbf{f,\beta)}=\sum_{j=1}^{m}\sum_{v=(t-k)\vee1}^{t\wedge T}%
\sum_{q=v}^{v+k}\beta_{j,q-v+1}\beta_{j,t-v+1}f_{q}.
\]
Then, calling $\mathbf{b(f,\beta)}=(b_{1}\mathbf{(f,\beta)},...,b_{T+k}%
\mathbf{(f,\beta)})^{\prime}$
\begin{equation}
\mathbf{b(f,\beta)=D(\beta)f,} \label{bf}%
\end{equation}
where $\mathbf{D}$ is given in (\ref{matD}). Then, from (\ref{abeta}) and
(\ref{bf}), \ equation \ (\ref{eqfun}) can be also written as $\sum_{j=1}%
^{m}\mathbf{C}_{j}(\alpha_{j})\mathbf{\beta}_{j}=\mathbf{D}(\mathbf{\beta)f.}$
Then (\ref{g1}) follows

\subsection{Proof of (\ref{matrixC})}

To prove (\ref{matrixC}) it is enough to show that we can find $\lambda$ such
that $\ $%

\begin{equation}
\lambda(1+a_{1}/a_{2})-1=(\lambda b/a_{2})^{2} \label{fdes}%
\end{equation}
and%
\begin{equation}
|\lambda|(1+a_{1}/a_{2})<1. \label{desin}%
\end{equation}
In this case \ (\ref{matrixC}) holds with
\begin{equation}
c=|\lambda|(1+a_{1}/a_{2}). \label{cc}%
\end{equation}
According to (\ref{fdes}) $\lambda$ should satisfy%
\begin{equation}
(b^{2}/a_{2}^{2})\lambda^{2}-(1+a_{1}/a_{2})\lambda+1=0. \label{eq2gr}%
\end{equation}
\ A necessary and sufficient condition for the existence of a real solution of
this equation is \ that $(1+a_{1}/a_{2})^{2}-4b^{2}/a_{2}^{2}\geq0$\ which is
\ equivalent to%
\begin{equation}
a_{1}+a_{2}\geq2|b|. \label{sw}%
\end{equation}
To prove this is enough
\[
\sum_{j=1}^{m}\beta_{j,0}^{2}+\sum_{j=1}^{m}\beta_{j,1}^{2}\geq2\sum_{j=1}%
^{m}\beta_{j0}\beta_{j,1}%
\]
which is always true. Solving (\ref{eq2gr}) we get that one of the roots is%
\begin{align*}
\lambda &  =\frac{a_{2}(a_{1}+a_{2})}{2b^{2}}-\frac{a_{2}^{2}}{2b^{2}}\left(
\frac{(a_{1}+a_{2})^{2}}{a_{2}^{2}}-\frac{4b^{2}}{a_{2}^{2}}\right)  ^{1/2}\\
&  =\frac{a_{2}^{2}}{2b^{2}}\left(  \left(  1+\frac{a_{1}}{a_{2}}\right)
-\left(  (1+\frac{a_{1}}{a_{2}})^{2}-4\frac{b^{2}}{a_{2}^{2}}\right)
^{1/2}\right)
\end{align*}
and therefore$|\lambda|<a_{2}(a_{1}+a_{2})/2b^{2}$ and using (\ref{cc})\ and
(\ref{sw}) we get $|c|=|\lambda|\left(  1+(a_{1}/a_{2}\right)  )<(a_{1}%
+a_{2})^{2}/(2b^{2})\leq1$, proving \ (\ref{desin}).

\subsection{Derivation of (\ref{impa5}) and (\ref{impa6}) \label{Wpr9}}

Differentiating (\ref{sj}) \ with respect to $f_{t}$ and using (\ref{fineq2})
we get \
\begin{equation}
\frac{\delta s_{j}(\beta,\alpha,\mathbf{f)}}{\partial f_{t\text{ }}}%
=\frac{-s\sum_{v=(t-k)\vee1}^{t\wedge T}w_{j,v}(z_{r}-\alpha)\beta
_{j,t-v+1}+s\sum_{v=(t-k)\vee1}^{t\wedge T}\sum_{i=0}^{k}w_{j,v}\beta
_{j,i+1}\beta_{j,t-v+1}f_{v+i}}{\sum_{h=(t-k)\vee1}^{t\wedge T}w_{j,h}%
r_{j,h}^{2}\ }. \label{deq1}%
\end{equation}

Differentiating (\ref{lossrob}) \ with respect to $f_{t}$ we get
\begin{equation}
\sum_{j=1}^{m}s_{j}\frac{\delta s_{j}(\beta,\alpha,\mathbf{f)}}{\partial
f_{t\text{ }}}=0, \label{deq2}%
\end{equation}
\bigskip and then, \ from (\ref{deq1}) and (\ref{deq2}) \ we get%

\[
\sum_{j=1}^{m}\sum_{v=(t-k)\vee1}^{t\wedge T}W_{j,t,v}(z_{r}-\alpha
)\beta_{j,t-v+1}=\sum_{j=1}^{m}\sum_{r=(t-k)\vee1}^{t\wedge T}\sum_{i=0}%
^{k}W_{j,t,v}\beta_{j,i+1}\beta_{j,t-v+1}f_{v+i},
\]
where $W_{j,t,v}$ is given by (\ref{WW}). This equation can also be written
as
\begin{equation}
a_{t}(\mathbf{f},\mathbf{\beta)=}b_{t}(\mathbf{f},\mathbf{\beta),}
\label{impa0}%
\end{equation}
where%
\[
a_{t}(\mathbf{f},\mathbf{\beta,\alpha)}=\ \sum_{j=1}^{m}\sum_{v=(t-k)\vee
1}^{t\wedge T}W_{j,t,v}(z_{r}-\alpha)\beta_{j,t-v+1}%
\]
and%
\[
b_{t}(\mathbf{f},\mathbf{\beta)=}\sum_{j=1}^{m}\sum_{r=(t-k)\vee1}^{t\wedge
T}\sum_{i=0}^{k}W_{j,t,v}\beta_{j,i+1}\beta_{j,t-v+1}f_{v+i}%
\]
Putting $q=t-v+1$ we get
\begin{align}
a_{t}(\mathbf{f},\mathbf{\beta,\alpha)}  &  =\ \sum_{j=1}^{m}\sum
_{q=1\vee(t-T+1)}^{(k+1)\wedge t}W_{j,t-q+1}(\mathbf{f,\beta}_{j}%
\mathbf{,\alpha},s)(z_{j,t-q+1}-\alpha_{j})\beta_{j,q}\nonumber\\
&  =\sum_{j=1}^{m}\mathbf{C}_{j}(\mathbf{f,\beta}_{j}\mathbf{,s}%
)\mathbf{\beta}_{j}\mathbf{,} \label{impa1}%
\end{align}
where $\mathbf{C}_{j}(\mathbf{f,\beta}_{j},\mathbf{s)}$ \ is the
$(T+k)\times(k+1)$ defined in (\ref{CC}). Putting $v+i=q$ we get%
\begin{align}
b_{t}(\mathbf{f},\mathbf{\beta,\alpha)}  &  =\sum_{j=1}^{m}\sum_{v=(t-k)\vee
1}^{t\wedge T}\sum_{q=v}^{v+k}W_{j,t,v}\beta_{j,q-v+1}\beta_{j,t-v+1}%
f_{q}\nonumber\\
&  =\mathbf{D}(\mathbf{f,\beta,}s)\mathbf{f}, \label{impa2}%
\end{align}
\ where $\mathbf{D}(\mathbf{f,\beta,\alpha,s})\ $is the $(T+k)\times(T+k)$
matrix \ defined in (\ref{DD}) \ and $\mathbf{s=}(s_{1},...s_{m}).$ Then from
(\ref{impa0}), (\ref{impa1}) and (\ref{impa2}) we derive (\ref{impa5}).
Differentiating \ (\ref{sj}) with respect to $\beta_{j,i}$ and $\alpha_{j},$
we get%

\[
\frac{1}{T}\sum_{t=1}^{T}\psi\left(  \frac{z_{j},_{v-}\sum_{i=0}^{k}%
\beta_{j,i+1}f_{v+i}-\alpha_{j}}{s_{j}}\right)  \left(  -s_{j}f_{v+i-1}%
-r_{j,v}\frac{\partial s_{j}}{\partial\beta_{j,i}}\right)  \ =0
\]
\[
\frac{1}{T}\sum_{t=1}^{T}\psi\left(  \frac{z_{j,v-}\sum_{i=0}^{k}\beta
_{j,i+1}f_{v+i}-\alpha_{j}}{s_{j}}\right)  \left(  -s_{j}f_{v+i-1}%
-r_{j,v}\frac{\partial s_{j}}{\partial\alpha_{j}}\right)  \ =0.
\]
Then putting $\partial s_{j}/\partial\beta_{j,i}=0,1\leq i\leq k+1$ and
$\partial s_{j}/\partial\alpha_{j}=0$ by (\ref{WW})\ we get the following
equations%
\begin{equation}
\sum_{t=1}^{T}\ w_{j,v}\left(  z_{j,v-}\sum_{i=0}^{k}\beta_{j,i+1}%
f_{v+i}-\alpha_{j}\right)  f_{v+i-1}=0,1\leq i\leq k+1\ \label{beteq}%
\end{equation}
and%
\begin{equation}
\sum_{t=1}^{T}w_{j,v}\left(  z_{j,v-}\sum_{i=0}^{k}\beta_{j,i+1}f_{v+i}%
-\alpha_{j}\right)  =0. \label{alfaeq}%
\end{equation}
From (\ref{beteq}) and (\ref{alfaeq}) equation (\ref{impa6}) follows immediately.

\bigskip\textbf{References}

Ahn, S. K. and Reinsel, G. C. (1988). Nested reduced-rank autoregressive
models for multiple time series, \textit{Journal of the American Statistical
Association}, \textbf{83}, 849--856.

Ahn, S. K. and Reinsel, G. C. (1990). Estimation for partially nonstationary
multivariate autoregressive models, \textit{Journal of the American
Statistical Association}, \textbf{85}, 813--823.

Box, G.E.P. and Tiao, G. C. (1977). A canonical analysis of multiple time
series, \textit{Biometrika}, \textbf{64}, 355--365.

Brillinger, D. R. (1981). \textit{Time Series Data Analysis and Theory},
Expanded edition, Holden-Day, San Francisco.

Davies, P.L. (1987), Asymptotic Behavior of S-Estimators of Multivariate
Location Parameters and Dispersion Matrices, The Annals of Statistics, 15, 1269-1292.

Forni, M., Hallin, M., Lippi, M. and Reichlin, L. (2000). The generalized
dynamic factor model: Identification and estimation, \textit{The Review of
Economic and Statistics}, \textbf{82}, 540--554.

Ku, W., R.H. Storer, and C. Georgakis (1995) Disturbance detection and
isolation by dynamic principal component analysis. \emph{Chemometrics and
Intelligent Laboratory Systems}, 30, 179-196.

Lam, C. and Yao, Q. (2012) Factor modeling for high dimensional time series:
Inference for the number of factors, \emph{The Annals of Statistics}, 40, 2, 694-726.

Maronna, R.A. (2005) Principal Components and Orthogonal Regression Based on
Robust Scales, \emph{Technometrics,} \ 47, 264-273.

Maronna, R. A., Martin, R. D., and Yohai, V. J. (2006) \textit{Robust
Statistics, }Wiley, Chichester.

Pe\~{n}a, D. and Box, G.E.P. (1987). Identifying a simplifying structure in
time series, \textit{Journal of the American Statistical Association},
\textbf{82}, 836--843.

Pe\~{n}a, D. and Poncela, P. (2006). Nonstationary dynamic factor analysis,
\textit{Journal of Statistical Planning and Inference}, 136,4, 1237-1256.

Reinsel, G. C. and Velu, R. P. (1998). \textit{Multivariate Reduced-Rank
Regression}, Springer, \ New York.$_{{}}$

Rousseeuw, P.J. and Yohai, V. (1984), \textquotedblleft Robust Regression by
Means of S estimators\textquotedblright, in Robust and Nonlinear Time Series
Analysis, edited by J. Franke, W. H\"{a}rdle, and R.D. Martin, Lecture Notes
in Statistics 26, Springer Verlag, New York, 256-274.

Seber, G.A.F. (1984) \textit{Multivariate observations}, Wiley, New York.

Shumway, R.H. \& Stoffer, D. S. (2000). \textit{Time Series Analysis and Its
Applications}. New York: Springer.

Stock, J.H. and Watson, M.W. (1988) Testing for Common Trends. \textit{Journal
of the American Statistical Association}, \textbf{83}, 1097-1107.

Tiao, G. C. and Tsay, R.S. (1989). Model specification in multivariate time
series. \textit{Journal of the Royal Statistical Society, Series B},
\textbf{51}, 157- 195.
\end{document}